\documentstyle{amsppt} 
\magnification=\magstep1
\NoRunningHeads
\NoBlackBoxes
\parindent=1em 
\vsize=7.4in
%macros

 %triple norm of #1

%
\topmatter
\title
Local Theory of Frames and Schauder Bases for Hilbert Space                
\endtitle
\author
Peter G. Casazza
\endauthor
\address
Department of Mathematics,
The University of Missouri,
Columbia, Missouri 65211
\endaddress
\email
pete\@casazza.math.missouri.edu
\endemail
\thanks
This research was supported by NSF DMS 9701234.  Part of this research was 
conducted while the author was a visitor at the "Workshop on Linear Analysis
and Probability", Texas A\&M University.
\endthanks
%\subjclass
%\endsubjclass
%
\abstract
We develope a local theory for frames on finite dimensional Hilbert spaces.
We show that for every frame $(f_{i})_{i=1}^{m}$ for an n-dimensional Hilbert
space, and for every $\epsilon > 0$, there is a subset $I\subset \{1,2,\ldots,m\}$ with $|I|\ge (1-\epsilon )n$ so that $(f_{i})_{i\in I}$ is a Riesz basis for its span with Riesz basis constant a function of $\epsilon$,  the frame bounds, and $(\|f_{i}\|)_{i=1}^{m}$, but independent of m and n. We also construct an example of a normalized frame for a Hilbert space $H$ which contains a subset which forms a Schauder basis for $H$, but contains no subset which is a Riesz basis for $H$.  We give examples to show that all of our results are best possible, and that all parameters are necessary.     
\endabstract
\endtopmatter
\document
\baselineskip=15pt
\heading{1.Introduction}
\endheading
\vskip10pt

Casazza and Christensen \cite{3} have shown that there is a tight frame for 
a Hilbert space which does not contain a Riesz basis.  Later, they observed \cite{4} that this frame does not even contain a subset which is a permutation
of a Schauder basis.  It follows from these results that there are normalized
frames for an n-dimensional Hilbert space $H_{n}$ (with quite good frame bounds of 1/2 and 2)  so that any subset of the
frame which forms a Riesz basis for $H_{n}$ has Riesz basis constant at least
$\sqrt{n}$.  That is, even a ``good'' frame for $H_{n}$ need not contain a subset 
which forms a ``good'' Riesz basis for $H_{n}$.  However, we will show that such
frames always contain a subset which is a good Riesz basis for a subspace of $H_{n}$ whose dimension is a percentage (arbitrarily close to one) of n.  We will produce similar results for Schauder bases for $H_{n}$, but now the Riesz basis 
constant will also depend upon the Hilbertian constant of the basis (and this is
a necessary constraint).   We also give examples to show that all our results are best possible and all the parameters are necessary.  Finally, we construct the first example of a normalized frame for a Hilbert space which contains a Schauder basis for $H$ but does not contain a Riesz basis for $H$.  This means
that our sequence is a normalized frame for $H$ and contains a subsequence
which is a Schauder basis for $H$, but any subset of the frame which is a
Schauder basis is no longer a frame (Since separated sets which are frames
are automatically Riesz bases for $H$ with frame bounds equal to the 
square of the Riesz
basis bounds). 

Our work relies heavily on some deep results of Bourgain and Tzafriri \cite{2} on 
restricted invertibility of linear operators acting on finite dimensional ${\ell}_{p}-$spaces.  For completeness, we will state the result from \cite{2} which is used in this paper.  We will denote by $(e_{i})_{i\in I}$ an 
orthonormal basis for a finite or infinite dimensional Hilbert space.

\proclaim{Theorem 1.1 (Bourgain and Tzafriri)}
There is a constant $c > 0$ so that, whenever $T:{{\ell}_{2}}^{n}\rightarrow {{\ell}_{2}}^{n}$ is a linear operator for which $\|Te_{i}\| = 1$; for all $1\le i \le n$, then there exists a subset $\sigma$ of $\{1,2,\ldots,n\}$ of cardinality $|\sigma| \ge \frac{cn}{\|T\|^{2}}$ so that
$$
\|\sum_{j\in \sigma} a_{j}Te_{j}\|^{2} \ge c\ \sum_{j\in \sigma}|a_{j}|^{2} 
$$
for any choice of scalars $(a_{j})_{j\in \sigma}$.
\endproclaim

This paper explores the relationship between frames and the local
theory of Banach spaces.  We now direct some comments towards the
reader interested in further explorations of these connections.  The result
of Bourgain-Tzafriri above fails for ${\ell}_{p}$ as stated (see the discussion at the end of section 2).  But,
with slightly stronger hypotheses, it can be done for
$1< p\le 2$ (see Theorem 7.2, \cite{2}).  Theorem 2.1 below
was certainly known to Bourgain-Tzafriri and to specialists in the
area, but does not seem to have been formally written down.  The
corresponding result of Theorem 2.1 (even with the stronger hypotheses
needed to get the Bourgain-Tzafriri result above) is 
unknown for  ${\ell}_{p}$, 
$p\not= 2$.  The problem is that to pass from having a ``fixed proportion'' 
of your set of vectors being well equivalent to the unit vectors in 
${\ell}_{p}$ to having an arbitrarily close to one proportion with this
property, requires being able to produce a good projection onto
your set of vectors.   In a Hilbert space, this property is for free, while
in ${\ell}_{p}$ there may not be such projections in general.  The 
arguements in this paper are similar to the so-called proportional 
Dvoretzky-Rogers factorization, as used for example by Szarek and
Talagrand \cite{9}.  The result from \cite{9} was improved by Giannopoulous \cite{5}.  Also, one can see these ideas in the paper of Bourgain and
Szarek \cite{1}.  Finally, in a paper in preparation, Litvak amd
Tomczak-Jaegermann \cite{6} describe the Dvoretzky-Rogers factorization
for non-symmetric bodies which give even stronger results than some of
ours, but are much more technical.  Finally, there are connections between
frames and convex geometry ralating tight frames to the 
so-called John's decomposition.  This is a bit technical for this paper and we refer the interested reader to \cite{7,10}.  

\heading{2. Local Theory of Schauder Bases}
\endheading
\vskip10pt

We say that two sets of vectors $(f_{i})_{i\in I}$ and $(g_{i})_{i\in I}$ are {\bf K-equivalent} if for every set of scalars $(a_{i})_{i\in I}$ we have,
$$
K^{-1}\|\sum_{i\in I}a_{i}f_{i}\| \le \|\sum_{i\in I}a_{i}g_{i}\|    \le K\|\sum_{i\in I}a_{i}f_{i}\|.
$$ 
A sequence $(f_{i},f_{i}^{*})_{i\in I}$ in $H$ is called a {\bf biorthogonal system with constant
d} if $\|f_{i}\| = 1$, and $\|f_{i}^{*}\| \le d^{-1}$, for all $i\in I$, and $<f_{i}^{*},f_{j}> = {\delta}_{ij}$, for all $i,\ j\in I$.  This is equivalent to $(f_{i})_{i\in I}$ being a set of vectors in $H$ satisfying: 
$$
d = \text{inf}_{j}\text{inf}\{\|f_{j} - f\| : f\in \text{span}(f_{i})_{i\not= j}\} > 0.
$$
A sequence of vectors $(f_{i})_{i\in I}$ is a {\bf Hilbertian sequence} with {\bf Hilbertian constant L} if
$$
\|\sum_{i\in I}a_{i}f_{i}\| \le L\left ( \sum_{i\in I}|a_{i}|^{2} \right )^{1/2},
$$
for all sequences of scalars $(a_{i})_{i\in I}$.  A sequence of vectors $(f_{i})_{i\in I}$ is a {\bf Besselian sequence} with {\bf Besselian constant B} if
$$
B \|\sum_{i\in I}a_{i}f_{i}\| \ge \left ( \sum_{i\in I}|a_{i}|^{2} \right )^{1/2},
$$
for all sequences of scalars $(a_{i})_{i\in I}$.  The sequence $(f_{i})_{i\in I}$ is called a {\bf Riesz basis} for its span with {\bf Riesz basis constant M} if
$$
\frac{1}{M}\left ( \sum_{i\in I} |a_{i}|^{2} \right )^{1/2} \le \|\sum_{i\in I}a_{i}f_{i}\| \le M \left ( \sum_{i\in I} |a_{i}|^{2} \right )^{1/2},
$$
for any choice of scalars $(a_{i})_{i\in I}$.

In several places in the paper we need certain conditional bases 
for a Hilbert space.  We will write down these bases now without
verifying their properties.  The proof can be found, for example, in
\cite{8}.

Let $H=L^{2}[-{\pi},\pi ]$ and let $0<a<1/2$.  Then
$$
f_{2n}(x) = |x|^{-a}e^{-inx},\ \ \ f_{2n+1}(x)= |x|^{-a}e^{inx},\ \ 
n=0,1,2,\ldots
$$
is a Besselian but non-Hilbertian bounded basis of $H$.  Also,
$$
g_{2n}(x) = |x|^{a}e^{inx},| | | g_{2n+1}(x) = |x|^{a}e^{-inx},\ \ n=0,1,2,\ldots
$$
is a Hilbertian but non-Besselian bounded basis of $H$.

We now show that finite separated, bounded, Hilbertian sequences have large subsets which are Riesz bases for their span.  We will give an explanation for the inner workings of this proof right afterwards.  

\proclaim{Theorem 2.1}
There is a function $g(x,y,z):R^{3} \rightarrow R^{+}$ with the following property:   
Let $(f_{i},f_{i}^{*})_{i=1}^{n}$ be any biorthogonal system with constant d,
 $0 < d \le 1$ and Hilbertian constant $L$ in an n-dimensional Hilbert space $H_{n}$ and let $0 < \epsilon < 1$.  Then there is a subset ${\sigma}\subset \{1,2,\ldots,n\}$, with $|\sigma| \ge (1-\epsilon )n$ so that $(f_{i})_{i\in \sigma}$ is a Riesz basis for its span with Riesz basis constant $g(\epsilon ,\ d,\ L)$. 
\endproclaim

\demo{Proof}
By defining an operator $T:{\ell}_{2}^{n} \rightarrow {\ell}_{2}^{n}$ by $T(e_{i}) = f_{i}$ and letting $b = \frac{cd^2}{L^{2}}$ in Theorem 1.1, we obtain a set ${\sigma}_{1}\subset \{1,2,\ldots,n\}$ with $|{\sigma}_{1}| = bn$ so that
$$
\|\sum_{j\in {\sigma}_{1}}a_{j}f_{j}\| \ge c \left ( \sum_{j\in {\sigma}_{1}}|a_{j}|^{2} \right )^{1/2},
$$
for any choice of scalars $(a_{j})_{j\in {\sigma}_{1}}$.  Here, and for the rest of this proof, to simplify notation we will ignore the fact that $bn$ may not actually be an integer.  By working with the greatest integer function, we can make this more exact, but the notation becomes unnecessarily cumbersome.  Let $P_{1}$ be the orthogonal projection of $H_{n}$ onto $\text{span}_{j\in {\sigma}_{1}}f_{j}$.  By the definition of the biorthogonal constant, we have
$$
\|(I - P_{1})f_{j}\| \ge d,  \tag 3.1
$$
for all $j\in {\sigma}_{1}^{c}$.  Define an operator 
$$
T_{1}:H_{|{\sigma}_{1}^{c}|}\rightarrow H_{|{\sigma}_{1}^{c}|}
$$
by
$$
T_{1}e_{i} = \frac{(I-P_{1})f_{i}}{\|(I-P_{1})f_{i}\|},
$$
 for all $i\in {\sigma}_{1}^{c}$, where $(e_{i})$ is an orthonormal basis
for $H_{|{\sigma}_{1}^{c}|}$.  Note that by (3.1) we have for all sequences
of scalars $(a_{i})$,
$$
\|\sum_{i\in {\sigma}_{1}^{c}} a_{i}T_{1}e_{i}\| \le 
\frac{1}{d} \|\sum_{i\in {\sigma}_{1}^{c}} a_{i}(I-P_{1})f_{i}\| \le
\frac{1}{d}\|\sum_{i\in {\sigma}_{1}^{c}}a_{i}f_{i}\| \le \frac{L}{d}.
$$
Thus, $\|T_{1}\|\le \frac{L}{d}$.    

Now we apply Theorem 1.1 to the operator $T_{1}$ and obtain  a set ${\sigma}_{2}\subset {\sigma}_{1}^{c}$ with $|{\sigma}_{2}| = b(1-b)n$ so that
$$
\|\sum_{j\in {\sigma}_{2}} a_{j}f_{j}\| \ge cd \left ( \sum_{j\in {\sigma}_{2}}|a_{j}|^{2} \right )^{1/2},
$$
for any choice of scalars $(a_{j})_{j\in {\sigma}_{2}}$.  Let $P_{2}$ be the orthogonal projection of $(I-P_{1})H_{n}$ onto $\text{span}_{j\in {\sigma}_{2}} f_{j}$ and observe that $(I-P_{2})(I-P_{1})$ is the orthogonal projection of $H_{n}$ onto the orthogonal complement of $\text{span}_{j\in {\sigma}_{1} \cup {\sigma}_{2}}f_{j}$, so again by the definition of the biorthogonal constant, we have
$$
\|(I-P_{2})(I-P_{1})f_{j}\| \ge d,
$$
for all $j\in ({\sigma}_{1} \cup {\sigma}_{2})^{c}$.  We continue to get disjoint sets $({\sigma}_{i})_{i=1}^{m}$  and orthogonal projections $(P_{i})_{i=1}^{m}$ satisfying:

(1)  $[1 - (1-b)^{m-1}] \ge 1 - \epsilon$;

(2)  $|{\sigma}_{i}| \ge b(1-b)^{i-1}n$, for all $1\le i \le m$;

(3)  $(I-P_{i})(I-P_{i-1})\cdots(I-P_{1})$ is the orthogonal projection of $H_{n}$ onto the orthogonal complement of $\text{span}\{f_{j}:j\in \cup_{k=1}^{i}{\sigma}_{k}\}$, for all $1\le i \le m$;

(4)  $\|(I-P_{i})f_{j}\| \ge d$, for all $j\in (\cup_{k=1}^{i}{\sigma}_{k})^{c}$, for all $1\le i \le m$;

(5)  For any $1\le i \le m$ and any choice of scalars $(a_{j})_{j\in {\sigma}_{i}}$ we have
$$
\|\sum_{j\in {\sigma}_{i}}a_{j}(I - P_{i})(I - P_{i-1})\cdots (I - P_{1})f_{j}\| \ge cd \left ( \sum_{j\in {\sigma}_{i}}|a_{j}|^{2} \right )^{1/2}.
$$
Now let $\sigma = \cup_{i=1}^{m}{\sigma}_{i}$ and note that
$$
|\sigma| = b\sum_{i=0}^{m-1}(1-b)^{i}n = [(1-(1-b)^{m-1}]n \ge (1-\epsilon)n.
$$
It remains to show that $(f_{i})_{i\in \sigma}$ is Besselian with constant a function of the stated parameters.  For later reference in the proof of
Theorem 4.3, the reader should note that the rest of the proof relies
only on that fact that we have a disjoint family of subsets of $\{1,2,
\ldots,n\}$ satisfying (1) through (5) above.  To see this, choose real numbers $r > 2,$ and $a$ so that $2L < cd(r - 1)$ and $r^{m+1}a < 1$, and choose any set of scalars $(a_{j})_{j\in \sigma}$ with 
$$
\sum_{j\in \sigma}|a_{j}|^{2} = 1.  \tag 3.2
$$
Now choose $1 \le i_{o} \le m$   
largest so that 
$$
\left ( \sum_{j\in {\sigma}_{i_{o}}}|a_{j}|^{2}\right )^{1/2}\ge r^{m-i_{o}}a.  \tag 3.3
$$
Such an $i_{o}$ must exist for otherwise,
$$
\left ( \sum_{j\in \sigma}|a_{j}|^{2} \right )^{1/2} \le \sum_{i=1}^{m}\left ( \sum_{j\in {\sigma}_{i}}|a_{j}|^{2}\right )^{1/2} \le 
$$
$$
\sum_{i=1}^{m}r^{i}a \le r^{m+1}a < 1, 
$$
contradicting (3.1) above.  Now we have,
$$
\|\sum_{j\in \sigma}a_{j}f_{j}\| \ge \|\sum_{i=1}^{i_{o}}\sum_{j\in {\sigma}_{i}}a_{j}f_{j}\| - \sum_{i=i_{o}+1}^{m}\|\sum_{j\in {\sigma}_{i}}a_{j}f_{j}\| \tag 3.4
$$
$$
\ge \|(I - P_{i_{o}-1})(I - P_{i_{o}-2})\cdots (I - P_{1})\left ( \sum_{i=1}^{i_{o}}\sum_{j\in {\sigma}_{i}}a_{j}f_{j}\right ) \| -     
\sum_{i=i_{o}+1}^{m}L\left ( \sum_{j \in {\sigma}_{i}}|a_{j}|^{2} \right ) ^{1/2}. 
$$
By our choice of $i_{o}$ largest satisfying (3.3) and from our construction, and (5) above we can continue inequality (3.4) as
$$
\ge \|\sum_{j\in {\sigma}_{i_{o}}}a_{j}(I - P_{i_{o}})(I - P_{i_{o}-1})\cdots (I - P_{1})f_{j}\| - L \sum_{i=i_{o}+1}^{m}r^{m-i}a 
$$
$$
\ge cd\left ( \sum_{j\in {\sigma}_{i_{o}}}|a_{j}|^{2}\right ) ^{1/2} - L \frac{r^{m-i_{o}}-1}{r-1} 
$$
$$
\ge cdr^{m-i_{o}}a - \frac{L}{r-1}r^{m-i_{o}}a \ge \frac{L}{r-1}r^{m-i_{o}}a,  
$$
where the last inequality follows from our choice of $r$.  Since $r^{m+1}a<1$, it follows that for every sequence of scalars $(a_{j})_{j\in \sigma}$ we have
$$
\|\sum_{j\in \sigma}a_{j}f_{j}\| \ge \frac{L}{r-1}a.
$$
This completes the proof of the Theorem.       
\enddemo

We feel that a discussion of the inner workings of the proof of Theorem 2.l is in order since on the surface such a proof should not work.  That is, we divided a
birothogonal system into subsets each of which is a good Riesz basis for its span and then took the union of these sets to get a larger Riesz basis.  Normally, such a process would fail for a biorthogonal system since our assumption is only that each
vector is far from the span of the others while we need that the span of certain subsets are far from the span of the others.  What is actually happening is the following.  We take an orthogonal projection $P$ onto the span of a subset $(f_{i})_{i\in \Delta}$ of our set of vectors and use 
``biorthogonality'' to discover that the vectors $((I - P)f_{i})_{i\in {\Delta}^{c}}$ are well bounded below in norm and hence have a subset $(I - P)f_{i})_{i\in {\Delta}_{1}}$ forming a good Riesz basis for their span.  Since $(f_{i})_{i\in {\Delta}_{1}}$ has Hilbertian constant $L$, it follows that $(f_{i})_{i\in {\Delta}_{1}}$ is also a good Riesz basis.  i.e.  $((I - P)f_{i})_{i\in {\Delta}_{1}^{c}}$ is well equivalent to $(f_{i})_{i\in {\Delta}_{1}}$.  It is not hard to see that this implies that the span of $(f_{i})_{i\in {\Delta}_{1}}$ is a ``good"
'' distance from the span of
$(f_{i})_{i\in {\Delta}}$, which is what we need.     
 
A sequence $(f_{i})_{i=1}^{m}$, with m finite or $m = \infty$, is a {\bf Schauder basis}
for $H_{m}$ if for every $f\in H_{m}$, there is a unique set of scalars $(a_{i})_{i=1}^{m}$ so that
$$
f = \sum_{i=1}^{m}a_{i}f_{i}.
$$
In the finite dimensional case, this is not particularly interesting since this is equivalent to the sequence being linearly independent.  What is important in this 
case is a quantative measure of the behavior of the basis.  The {\bf basis constant K} of the Schauder basis $(f_{i})_{i=1}^{m}$ is the smallest constant satisfying:
$$
\|\sum_{i=1}^{n}a_{i}f_{i}\| \le K\|\sum_{i=1}^{m}a_{i}f_{i}\|,
$$
for every natural number $n < m$ and every choice of scalars $(a_{i})_{i=1}^{m}$.
It is easily checked that if $(f_{i})_{i\in I}$ is a Schauder basis with basis constant $K$, then $(f_{i})_{i\in I}$ is a separated set with constant $\ge \frac{1}{2K}$. To get a separated set of vectors which is not a Schauder basis, take any conditional Schauder basis $(f_{i})_{i=1}^{\infty}$ for $H$ and choose a permutation $\sigma$ of the natural numbers so that $(f_{{\sigma}(i)})_{i=1}^{\infty}$ is not a Schauder basis for $H$.  Then this set is still separated but is no longer a Schauder basis for $H$.  The next result is an immediate consequence of Theorem 2.1.
 
\proclaim{Corollary 2.2}
There is a function $g(x,y,z):R^{3} \rightarrow R^{+}$ with the following property:   
Let $(f_{i})_{i=1}^{n}$ be any normalized Schauder basis for an n-dimensional Hilbert space $H_{n}$ with basis constant $K$ and Hilbertian constant $L$, and let $0 < \epsilon < 1$.  Then there is a subset $\sigma \subset \{1,2,\ldots,n\}$, with $|\sigma| \ge (1-\epsilon )n$ so that $(f_{i})_{i\in \sigma}$ is a Riesz basis for its span with Riesz basis constant $g(\epsilon ,\ K,\ L)$.    
\endproclaim

Corollary 2.2 (and even Theorem 1.1) do not generalize to ${\ell}_{p}$.  For example, if $1\le p < 2$, there is a constant $K>0$ so that for all $n$,  ${\ell}_{p}^{2n}$ contains a normalized sequence $(x_{i})_{i=1}^{n}$ which is K-equivalent to the unit vector basis of ${\ell}_{2}^{n}$.  Defining $T:{\ell}_{p}^{2n} \rightarrow {\ell}_{p}^{2n}$ by $Te_{i} = x_{i}$, for $1\le i \le n$ and $Te_{i} = 0$, for $n+1 \le i \le 2n$, it follows that $\|T\| = 1$, yet we do not have a large
subset of $(x_{i})_{i=1}^{2n}$ which is well equivalent to the unit vector
basis of ${\ell}_{p}^{m}$.  This shows that our results do not work in general outside of Hilbert space.

It is well known that there are conditional Schauder bases (even Hilbertian ones) $(f_{i})_{i=1}^{\infty}$ for a Hilbert space.  This means that the Riesz basis constant of $(f_{i})_{i=1}^{n}$ goes to infinity with n.  The main point of Theorem 2.1 is that the Riesz basis constant is independent of the dimension of the Hilbert space.  That is, although $(f_{i})_{i=1}^{n}$ itself need not be a Riesz basis for $H_{n}$ with Riesz basis constant independent of n, at least it has a subset spanning a (arbitrarily close to one) percentage of the dimension of the space which is a Riesz basis for its span with Riesz basis constant independent of n (but of course a function of the percentage).  Each of the variables in the function $g(x,\ y, \ z)$ are necessary for Theorem 2.1 and Corollary 2.2 to hold. The preceeding discussion shows that $\epsilon$ is necessary in these results.  Bourgain and Tzafriri \cite{2} (the Remark on p-165) give an example of a Besselian Schauder basis $(f_{i})_{i=1}^{\infty}$ for $H$ which has no subset of positive density which is a Riesz basis for its closed linear span.  This means that $(f_{i})_{i=1}^{n}$ does not contain a percentage which is a Riesz basis with Riesz basis constant independent of n.  Finally, the separation assumption in Theorem 2.1 is necessary since otherwise we could consider $(e_{i}, \ e_{i})_{i=1}^{n}$ in $H_{2n}$ and have no subset at all which is a Riesz basis for more than half of $H_{2n}$.  If we want our set to be linearly independent, we can use $(e_{2i-1}, \ e_{2i-1} + \frac{1}{n}e_{2i})_{i=1}^{n}$ in $H_{2n}$ and easily observe that this is a linearly independent set spanning $H_{2n}$ for which any subset containing more than half the elements has Riesz basis constant $\sqrt{2}n$.

\heading{3. Every Frame is Equivalent to a Tight Frame}
\endheading
\vskip10pt

The results of this section have been part of the folklore in this area for some time, but do not seem to be broadly known.  Recall that a sequence $(f_{i})_{i\in I}$ in a Hilbert space $H$ is a {\bf frame} for $H$ with frame bounds $A,\ B$ if
$$
A\|f\|^{2} \le \sum_{i\in I}|<f,f_{i}>|^{2} \le B\|f\|^{2}, \ \ \forall f\in H.
$$
If $A = B$, we call this a {\bf tight frame}.  If $(f_{i})_{i\in I}$ is a frame, then defining $Sf = \sum_{i\in I}<f,f_{i}>f_{i}$, for all $f\in H$, we obtain an isomorphism of $H$ onto $H$.  $S$ is called the {\bf frame operator} for the frame.  This leads to the {\bf frame decomposition},
$$
f = SS^{-1}f = \sum_{i\in I}<f,S^{-1}f_{i}>f_{i} = \sum_{i\in I}<S^{-1}f,f_{i}>f_{i}, \ \ \forall \  f\in H.
$$
It follows that
$$
<f,S^{-1}f> = \sum_{i\in I}|<S^{-1}f,f_{i}>|^{2}.  \tag 3.1
$$
As a consequence of (3.1), we can see that a frame is tight if and only if the frame operator is a multiple of the identity. 
The frame operator $S$ is easily seen to be a positive operator on $H$ and therefore real powers of $S$ make good sense.  This leads to the following general result.   

\proclaim{Theorem 3.1}
Let $(f_{i})_{i\in I}$ be a frame for $H$ with frame operator $S$.  Then for
any real number $a$, $(S^{\frac{a-1}{2}}f_{i})_{i\in I}$ is also a frame for $H$ with frame operator $S^{a}$.
\endproclaim

\demo{Proof}
Since $S$ is a positive operator and an isomorphism of $H$ onto $H$, so is $S^{b}$ for any real number b.  Hence, $(S^{b}f_{i})_{i\in I}$ is a frame for $H$.  Letting $b = \frac{a-1}{2}$ we compute for all $f\in H$,
$$
\sum_{i\in I}<f,S^{b}f_{i}>S^{b}f_{i} = S^{b} \left ( \sum_{i\in I}<f,S^{b}f_{i}>f_{i} \right ) 
$$
$$
= S^{b}\left ( \sum_{i\in I}<S^{b}f,f_{i}>f_{i}\right ) = S^{b}S(S^{b}f) = S^{1+2b}f = S^{a}f.
$$
This shows that $(S^{\frac{a-1}{2}}f_{i})_{i\in I}$ is a frame for $H$ with frame operator $S^{a}$. 
\enddemo
Letting $\frac{a-1}{2} = -\frac{1}{2}$, we get that $a=0$ in Theorem 3.1.  This yields,

\proclaim{Corollary 3.2}
If $(f_{i})_{i\in I}$ is a frame with frame operator $S$, then $(S^{-1/2}f_{i})_{i\in I}$ is a frame with the identity as frame operator.  That is, for every $f\in H$,
$$
f = \sum_{i\in I}<f,S^{-1/2}f_{i}>S^{-1/2}f_{i}.
$$
Therefore, every frame is equivalent to a tight frame.
\endproclaim

There are many places in the literature on frames where authors find (or the reader is asked to find) ``tight frame'' examples of an existing example in frame theory.  Corollary 3.2 renders all this as unnecessary, despite its relatively soft proof.

\proclaim{Corollary 3.3}
A frame $(f_{i})_{i\in I}$ is a Riesz basis for $H$ if and only if $(S^{-1/2}f_{i})_{i\in I}$ is an orthonormal basis for $H$.
\endproclaim

\demo{Proof}
$(S^{-1/2}f_{i})_{i\in I}$ is an orthonormal bais for $H$ if and only if
$$
{\delta}_{i,j} = <S^{-1/2}f_{i},S^{-1/2}f_{j}> = <S^{-1}f_{i},f_{j}>.
$$
That is, $(S^{-1/2}f_{i})_{i\in I}$ is an orthonormal basis for $H$ if and only if $(S^{-1}f_{i},f_{i})_{i\in I}$ is a biorthogonal sequence in $H$.  But, it is well known \cite{11} that this is equivalent to $(f_{i})_{i\in I}$ being a Riesz basis for $H$.
\enddemo

\heading{4.  Local theory of Frames}
\endheading
\vskip10pt

Casazza and Christensen \cite{3,4} (also see Lemma 5.1 below) have shown that there exist tight frames $(f_{i})_{i=1}^{n+1}$ for $H_{n}$ with $1/2\le \|f_{i}\|\le 2$ for which any subset which spans $H_{n}$ has Riesz basis constant $\ge \frac{\sqrt{n-1}}{4}$.  That is, a frame for a finite dimensional Hilbert space (even a tight frame with good bounds on the norms of the frame elements) need not contain a subset which is a Riesz basis for the space with Riesz basis constant independent of the dimension of the space.  However, in this section we will show that such frames contain ``good'' Riesz bases for a subspace ``almost'' equal to the whole space.  These results are just an application of the results of Section 2.  To apply the results of Section 2, we need two elementary observations.  The first result relates the dimension of the space to the lower frame bound and and the maximum of the norms of the frame elements.  

\proclaim{Lemma 4.1}
Let $(f_{i})_{i\in I}$ be a frame for $H_{n}$ with lower frame bound $A$ and $\|f_{i}\| \le \delta$, for all $i\in I$.  Then
$$
n \le \frac{{\delta}^{2}}{A}|I|.
$$
\endproclaim

\demo{Proof}
For any $1\le j \le n$,
$$
A \le \sum_{i\in I}|<e_{j},f_{i}>|^{2}.
$$
Therefore,
$$
nA \le \sum_{j=1}^{n}\sum_{i\in I}|<e_{j},f_{i}>|^{2} =
$$
$$
\sum_{i\in I}\sum_{j=1}^{n}|<e_{j},f_{i}>|^{2} = \sum_{i\in I}\|f_{i}\|^{2} \le {\delta}^{2}|I|.
$$
\enddemo

Our next preliminary result relates the cardinality of the number of frame elements to the upper frame bound, the dimension of the space and the minimum of the norms of the frame elements.

\proclaim{Lemma 4.2}
Let $(f_{i})_{i\in I}$ be a frame for $H_{n}$ with upper frame bound $B$ and
$\alpha \le \|f_{i}\|$, for all $i\in I$.  Then
$$
|I| \le \frac{B}{{\alpha}^{2}}n.
$$
\endproclaim

\demo{Proof}
We compute,
$$
{{\alpha}^{2}}|I| \le \sum_{i\in I}\|f_{i}\|^{2} = \sum_{i\in I}\sum_{j=1}^{n}|<e_{j},f_{i}>|^{2} =
$$
$$
\sum_{j=1}^{n}\sum_{i\in I}|<e_{j},f_{i}>|^{2} \le \sum_{j=1}^{n}B\|e_{j}\|^{2} = nB.
$$
\enddemo

Now we are ready for the main result of this section.

\proclaim{Theorem 4.3}
There is a function $g(v,w,x,y,z):R^{5} \rightarrow R^{+}$ with the following property:   
Let $(f_{i})_{i=1}^{k}$ be any frame for an n-dimensional Hilbert space $H_{n}$ with frame bounds $A,\ B$, ${\alpha} \le \|f_{i}\| \le {\beta}$, for all $1\le i \le k$, and let $0 < \epsilon < 1$.  Then there is a subset $\sigma \subset \{1,2,\ldots,k\}$, with $|\sigma| \ge (1-\epsilon )n$ so that $(f_{i})_{i\in \sigma}$ is a Riesz basis for its span with Riesz basis constant $g(\epsilon ,\ A,\ B,\ {\alpha},\ {\beta})$.    
\endproclaim

\demo{Proof}
By Lemma 4.2,
$$
k \le \frac{B}{{\alpha}^{2}}n.
$$
Now choose ${\delta} > 0$, a function of our stated perameters, so that
$$
\frac{{\delta}^{2}}{A} \frac{B}{{\alpha}^{2}} \le \frac{\epsilon}{2}. \tag 4.1
$$
Since a frame is Hilbertian with constant $\le B$,
by Theorem 1.1 there is a universal constant $c$ and a constant $d = c/B^{2}$ so that we can choose ${\sigma}_{1} \subset \{1,2,\ldots,k\}$ with $|{\sigma}_{1}| \ge dn$ and
$$
\|\sum_{i\in {\sigma}_{1}}a_{i}f_{i}\| \ge c \left ( \sum_{i\in {\sigma}_{1}}|a_{i}|^{2} \right )^{1/2},
$$
for all choices of scalars $(a_{i})_{i\in {\sigma}_{1}}$.  Let $P_{1}$ be the orthogonal projection of $H_{n}$ onto the span of $(f_{i})_{i\in {\sigma}_{1}}$.  If
$$
|\{i\in {\sigma}_{1}^{c} : \|(I - P_{1})f_{i}\| \ge \delta \}| \ge \frac{\epsilon}{2}n,
$$
then applying Theorem 1.1 again we can find ${\sigma}_{2} \subset {\sigma}_{1}^{c}$ with 
$|{\sigma}_{2}| \ge \frac{d}{{\delta}^{2}}\frac{\epsilon}{2}n$, so that

$$
\|\sum_{i\in {\sigma}_{2}}a_{i}f_{i}\| \ge c{\delta} \left ( \sum_{i\in {\sigma}_{2}}|a_{i}|^{2} \right ) ^{1/2},
$$
for all choices of scalars $(a_{i})_{i\in {\sigma}_{2}}$.  Let $P_{2}$ be the orthogonal projection of $H_{n}$ onto the span of $(f_{i})_{i\in {\sigma}_{2}}$, and check if
$$
|\{i \in ({\sigma}_{1} \cup {\sigma}_{2})^{c} : \|(I-P_{2})(I-P_{1})f_{i}\| \ge \delta \}| \ge \frac{\epsilon}{2}n.
$$
We continue this construction stopping it after m steps as soon as one of the following holds:

(1)  Letting
$$
{\sigma}_{m+1} = \{ i\in \left ( \cup_{j=1}^{m}{\sigma}_{j}\right )^{c} : \|(I-P_{m})(I-P_{m-1})\cdots(I-P_{1})f_{i}\| \ge \delta \}, \tag 4.2
$$
then $m$ is the first natural number so that:
$$
|{\sigma}_{m+1}| \le \frac{\epsilon}{2}n,
$$
or 

(2)  $dn + (m-1)\frac{d}{{\delta}^{2}}\frac{\epsilon}{2}n \ge (1-\epsilon)n$.

Now, let 
$$
\sigma = \cup_{j=1}^{m}{\sigma}_{j} .
$$
We finish the proof in two steps.

\proclaim{Step I}
$|\sigma| \ge (1-{\epsilon})n$.
\endproclaim

There are two cases to be examined here.

\proclaim{Case I}
$dn + m\frac{d}{{\delta}^{2}}\frac{\epsilon}{2}n \ge (1-\epsilon)n$
\endproclaim

In this case,
$$
|\sigma | = \sum_{j=1}^{m}|{\sigma}_{j}| = dn + (m-1)\frac{d}{{\delta}^{2}}\frac{\epsilon}{2}n \ge (1 - {\epsilon})n.
$$

\proclaim{Case II}
$m$ is the first natural number so that:
$$
|\{ i\in \left ( \cup_{j=1}^{m}{\sigma}_{j}\right )^{c} : \|(I-P_{m})(I-P_{m-1})\cdots(I-P_{1})f_{i}\| \ge \delta \}| \le \frac{\epsilon}{2}n,
$$
\endproclaim   

In this case, let $P_{m+1}$ be the orthogonal projection of $H_{n}$ onto the span of $(f_{i})_{i\in {\sigma}_{m+1}}$. It follows that,
$$
\|(I-P_{m+1})(I-P_{m})(I-P_{m-1})\cdots(I-P_{1})f_{i}\| \le \delta,
$$
for all $i\in ({\sigma} \cup {\sigma}_{m+1})^{c}$.  Now applying Lemma 4.1 and then Lemma 4.2 and then inequality (4.1) we have,
$$
\text{dim}\left ( \text{span}(f_{i})_{i\in ({\sigma} \cup {\sigma}_{m+1})^{c}}\right ) \le \frac{{\delta}^{2}}{A}k \le
\frac{{\delta}^{2}}{A} \frac{B}{{\alpha}^{2}}n \le \frac{\epsilon}{2}n.
$$
Combining this with inequality (4.2) yields
$$
\text{dim}\left ( \text{span}(f_{i})_{i\in {\sigma}^{c}}\right ) \le {\epsilon}n.
$$
Therefore, since $(f_{i})_{i\in I}$ spans $H_{n}$, it follows that $|\sigma|\ge (1-{\epsilon})n.$

The proof will be finished if we prove,
\proclaim{Step II}
$(f_{i})_{i\in \sigma}$ is a Riesz basis for its span with constant $g(v,w,x,y,z)$.
\endproclaim  
But, the (end of the) proof of Theorem 2.1 works here to show that our set is a good Riesz basis.  That is, the $({\sigma}_{i})$ above satisfy (3) through (5)
of the proof of Theorem 2.1, and hence from that proof for a good Riesz basis
for their span.  This completes the proof of the theorem.   
\enddemo

Again, the important point in Theorem 4.3 is that the Riesz basis constant is a function of the frame bounds, the max and min of the norms of the frame elements, and $\epsilon$, but is independent of the dimension of the space.  It is easily seen that all the parameters are necessary in Theorem 4.3.  Our earlier examples with Hilbertian Schauder bases show all this except the boundedness assumption.  But the frame given at the beginning of Section 5 below shows that the boundedness assumption is also necessary in theorem 4.3.

\heading{5.  Frames Containing Schuader Bases but not Riesz Bases}
\endheading
\vskip10pt

It is easy to construct a tight frame for a Hilbert space which contains a Schauder basis but does not contain a Riesz basis.  Just consider
$$
\{e_{1}, \frac{1}{\sqrt{2}}e_{2},\frac{1}{\sqrt{2}}e_{2},\frac{1}{\sqrt{3}}e_{3},\frac{1}{\sqrt{3}}e_{3},\frac{1}{\sqrt{3}}e_{3},\ldots\}.
$$
This frame has a subset $(\frac{1}{\sqrt{n}}e_{n})_{n=1}^{\infty}$ which is a Schauder basis for $H$.  But, any spanning subset of this frame is not bounded below in norm and hence is not a Riesz basis for $H$.  However, to construct an example of this type which is normalized is much more difficult, and has been open for quite a time.  We will give such an example below.  But we will first state the results needed for the example.  The first is due to Casazza and Christensen \cite{4}, Lemma 3 (This is not exactly what their lemma states, but it is what their proof yields).  

\proclaim{Lemma 5.1}
Let $(e_{i})_{i=1}^{n}$ be an orthonormal basis for an n-dimensional Hilbert space $H_{n}$.  Define
$$
f_{i} = e_{i}-\frac{1}{n}\sum_{j=1}^{n}e_{j}, \ \text{for all} \ \ i = 1,2,\ldots,n,
$$
and let
$$
f_{n+1} = \frac{1}{\sqrt{n}}\sum_{j=1}^{n}e_{j}.
$$
Then $(f_{j})_{j=1}^{n+1}$ is a frame for $H_{n}$ with bounds $A = B = 1$, and any subset of the frame which contains n-elements has basis constant greater than or equal to $\frac{\sqrt{n-2}}{4}$. 
\endproclaim

 We also need a particular example of a conditional Schauder basis for finite dimensional Hilbert spaces.
 
\proclaim{Lemma 5.2}
 There are universal constants $c,\ L$ so that for every $\epsilon > 0$ and every natural number $k$, there is a natural number $n$ and a normalized Hilbertian Schauder basis
$(g_{i})_{i=1}^{n}$ for $H_{n}$ with basis constant $c$ and Hilbertian constant $L$, and there is a subspace $E\subset H_{n}$ with dim E = k, and
$$
\sum_{i=1}^{n}|<g_{i},f>|^{2} \le {\epsilon}\|f\|^{2}, \ \ \forall \ f\in E.
$$
\endproclaim

\demo{Proof}
Let $(h_{i})_{i=1}^{\infty}$ be a normalized conditional Schauder basis for $H$ with basis constant $c$ and Hilbertian constant $L$.  Since $(h_{i})_{i = 1}^{\infty}$ is not a Riesz basis, it follows that for every $\epsilon > 0$ and every natural number $k$, there is a natural number $m$ and an vector $h\in \text{span}_{1\le i \le m}h_{i}$ with $\|h\| = 1$ and
$$
\sum_{i=1}^{m}|<h,h_{i}>|^{2} \le \frac{\epsilon}{k}.
$$
Let
$$
H_{n} = \left ( \sum_{j=1}^{k}\oplus H_{m} \right ) _{{\ell}_{2}}.
$$
Then the following sequence $(f_{ij})_{i=1,j=1}^{\ \ m \ \ \ k}$ forms a Schauder basis for $H_{n}$ with basis constant $c$ and Hilbertian constant $L$:
$$
\{(h_{1},0,\ldots,0),(h_{2},0,\ldots,0),\cdots,(h_{m},0,\ldots,0),
$$
$$
(0,h_{1},0,\ldots,0),\cdots,(0,h_{m},0,\ldots,0),\cdots,
$$
$$
(0,\ldots,0,h_{1}),\cdots,(0,\ldots,0,h_{m})\}.
$$
Let $f_{1} = (h,0,\ldots,0),f_{2} = (0,h,0,\ldots,0),\cdots,f_{k} = (0,\ldots,0,h)$, and let $E = \text{span}_{1\le i \le k}f_{i}$.  For any 
$$
f = \sum_{j=1}^{k}a_{j}f_{j} \in E,
$$
with
$$
\|f\|^{2} = \sum_{j=1}^{k}|a_{j}|^{2} = 1,
$$
we have that
$$
\sum_{i=1}^{m}\sum_{j=1}^{k}|<f,f_{ij}>|^{2} \le k\frac{\epsilon}{k} = \epsilon.
$$
This completes the proof of the lemma.
\enddemo

Now we are ready for the promised example.

\proclaim{Proposition 5.3}
There is a normalized frame $(f_{i})_{i=1}^{\infty}$ for $H$ which contains a Schauder basis but does not contain a Riesz basis (and hence does not contain any subset which is an unconditional basis for $H$).
\endproclaim

\demo{Proof}
Let ${\epsilon}_{m}$ decrease to 0, and by Lemma 5.2, choose $n_{m}$ and $E_{m}\subset H_{n_{m}}$, with $\text{dim}E_{m} = m$, and $(g_{i}^{m})_{i=1}^{n_{m}}$ satisfying Lemma 5.2.  Let
$$
H = \left ( \sum_{m=1}^{\infty}\oplus H_{n_{m}} \right )_{{\ell}_{2}}
$$
and let $P_{m}$ be the orthogonal projection of $H$ onto $E_{m}$.  Choose the tight frame $(f_{i}^{m})_{i=1}^{m+1}$ for $E_{m}\subset H$ as in Lemma 5.1 and let $(e_{i}^{m})_{i=1}^{n_{m}-m}$ be an orthonormal basis for $(I - P)H_{n_{m}}\subset H$.  Now, $\{(e_{i}^{m})_{i=1}^{n_{m}},\ (f_{i}^{m})_{i=1}^{m+1}\}$ forms a tight frame for $H_{n_{m}}\subset H$ and $(g_{i}^{m})_{i=1}^{n_{m}}$ is a Schauder basis for $H_{n_{m}}\subset H$ with basis constant $c$ and Hilbertian constant $L$.  Although not all the vectors here are normalized, they all have norms between 1 and 2, and so if we establish they satisfy the requirements of the Proposition, then normalizing them later will be sufficient.  From our observations above, the set of vectors 
$$
\{(g_{i}^{m})_{i=1}^{n_{m}},(e_{i}^{m})_{i=1}^{n_{m}},\ (f_{i}^{m})_{i=1}^{m+1}\}
$$
forms a frame in $H$ which contains the Schauder basis $\{(g_{i}^{m})_{i=1}^{n_{m}}\}_{m=1}^{\infty}$.  Now we need to show that this frame does not contain a Riesz basis for $H$.  Suppose that $(h_{i})_{i=1}^{\infty}$ is a subset of this frame which spans $H$.  For each $m = 1,2,\ldots$, let $(f_{i}^{m})_{i\in {\Delta}_{m}}$ be the elements from $(f_{i}^{m})_{i=1}^{m+1}$ contained in $(h_{i})_{i=1}^{\infty}$.  Then for each $m = 1,2,\ldots$ we have two possibilities.

\proclaim{Case I}
$|{\Delta}_{m}| = m$.
\endproclaim

In this case, by Lemma 5.1 we have that the Riesz basis constant of $(h_{i})_{i=1}^{\infty}$ is $\ge \frac{\sqrt{m-2}}{4}$.

\proclaim{Case II}
$|{\Delta}_{m}| < m$.
\endproclaim

In this case, $(f_{i}^{m})_{i\in {\Delta}_{m}}$ does not span $E_{m}$, so there is  a vector $g\in E_{m}$ with $\|g\| = 1$ and so that
$$
\sum_{i\in {\Delta}_{m}}|<g,f_{i}^{m}>|^{2} = 0.
$$
But, $g \perp \text{span}(e_{i}^{m})$ so
$$
\sum_{i\in {\Delta}_{m}}|<g,e_{i}^{m}>|^{2} = 0.
$$
Finally, by our construction,
$$
\sum_{i=1}^{n_{m}}|<g,g_{i}^{m}>|^{2} \le {\epsilon}_{m}.
$$
That is, in Case II, the Riesz basis constant of $(h_{i})_{i=1}^{\infty}$ is $\ge \frac{1}{{\epsilon}_{m}}$.

Combining Cases I and II for every $m$, we see that $(h_{i})_{i=1}^{\infty}$ is not a Riesz basis for $H$.      
\enddemo

\enddocument